\AddToHook{begindocument/before}{} 
\DeclareHookRule{begindocument}{hyperref}{before}{my-mdpi} 

\documentclass[preprints,article,accept,moreauthors,pdftex]{my-mdpi}

\firstpage{1} 
\pubvolume{11}
\issuenum{19}
\articlenumber{4218}
\pubyear{2023}
\copyrightyear{2023}
\externaleditor{Academic Editor: Andrea Scozzari}
\datereceived{12 September 2023} 
\daterevised{7 October 2023} 
\dateaccepted{8 October 2023} 
\datepublished{9 October 2023} 
\hreflink{https://\linebreak doi.org/10.3390/math11194218} 
\doinum{10.3390/math11194218}

\pdfoutput=1 

\usepackage{mathrsfs}
\usepackage{enumerate}

\Title{Pontryagin Maximum Principle for Incommensurate Fractional-Orders Optimal Control Problems}

\TitleCitation{Pontryagin Maximum Principle for Incommensurate Fractional-Orders Optimal Control Problems}

\Author{Fa\"{\i}\c{c}al Nda\"{\i}rou \orcidA{} 
and Delfim F. M. Torres *\orcidB{}}

\AuthorNames{Fa\"{\i}\c{c}al Nda\"{\i}rou and Delfim F. M. Torres}

\AuthorCitation{Nda\"{\i}rou, F.; Torres, D.F.M.}

\address[1]{Center for Research and Development in Mathematics and Applications (CIDMA),
Department of Mathematics, University of Aveiro, 3810-193 Aveiro, Portugal; 
faical@ua.pt}

\corres{\hangafter=1 \hangindent=1.05em \hspace{-0.82em}Correspondence: delfim@ua.pt; +351-234-370-668}

\abstract{We introduce a new optimal control problem where the controlled dynamical 
system depends on multi-order (incommensurate) fractional differential equations. 
The cost functional to be maximized is of Bolza type and depends on incommensurate Caputo
fractional-orders derivatives. We establish continuity and differentiability of the state solutions
with respect to perturbed trajectories. Then, we state and prove a Pontryagin maximum principle 
for incommensurate Caputo fractional optimal control problems. Finally, we give an example,
illustrating the applicability of our Pontryagin maximum principle.}



\begin{document}

\section{Introduction}

The celebrated Pontryagin maximum principle 
was first formulated in 1956, being widely regarded 
as the central result of optimal control theory \cite{MR0186436}.
The significance of the maximum principle lies in the fact 
that rather than maximizing over a function space, the problem 
is converted to a pointwise optimization. Various generalizations 
of the Pontryagin maximum principle are still being investigated
nowadays \cite{MR1360918,MR2536634,MR4626944}.

Fractional optimal control theory, as a branch of applied mathematics, 
deals with optimization issues for controlled fractional dynamics coupled 
with a performance index functional \cite{ButkovskiiI,ButkovskiiII}.
The development of the theory started shortly at the beginning of \mbox{21th 
century} with some pioneering works \cite{agarwal,jeli}, where necessary 
optimality conditions are derived using techniques from variational analysis. 
A major contribution from those works is weak versions of the Pontryagin 
maximum principle that consists in reducing the fractional optimal control 
problem (FOCP) to a fractional boundary value problem together with an 
optimality condition. Later, the theory emerged with the work of \cite{k}, 
by establishing a strong maximum principle of Pontryagin type, 
which introduces a Hamiltonian maximality condition 
in the set of necessary optimality conditions instead 
of the optimality condition from weak versions. Indeed, since 
this seminal work of 2014, the subject of maximum principle for 
FOCPs gained more interest. For instance, a new approach to the 
Pontryagin maximum principle for a problem involving a Lagrangian 
cost of Riemann--Liouville fractional integral subject to Caputo 
fractional dynamics is considered in \cite{lobo}. The authors of 
\cite{bourdin}, derived a Pontryagin maximum principle for FOCP 
defined by a general Bolza cost of fractional integral type 
with terminal constraints on the Caputo fractional dynamics. 
More recently, in 2022, Kamocki investigated an optimal 
control problem of multi-order fractional systems under 
a Lagrange-type functional, proving two existence results \cite{MR4449935}. 
In contrast with Kamocki, here we are interested in proving the necessary 
optimality conditions of Pontryagin type.

Another important development in control theory is the so-called Bellman's 
dynamic programming principle, which gives a necessary and sufficient condition 
of optimality \cite{Bardi}. In this direction, the approach for solving 
an optimal control problem can be achieved by determining a certain value 
function that might happen to be a viscosity solution to a Hamilton--Jacobi--Bellman (HJB) 
equation. The dynamic programming principle has been extended for fractional discrete-time 
systems by the authors of references \cite{czy,dzi}. Furthermore, an attempt to derive 
a fractional version of the Hamilton--Jacobi--Bellman (HJB) equation has been investigated 
in \cite{delfimm}. \textls[-25]{In \cite{mika1}, Gomoyunov studies an extension of the dynamic 
programming principle to the case of fractional-order dynamical systems. A new approach 
is required for such problems so that for every intermediate time $t \in (0, T)$, 
it is necessary to introduce an auxiliary optimal control problem (sub-problem) 
with this time $t$ considered as the initial one. Also, a derivation of the fractional 
version of the HJB equation is studied deeply \cite{mika1}. To sum up, it is important 
to mention that the Pontryagin maximum principle, along with Bellman's dynamic 
programming principle, is one of the most effective and fundamental tools 
for investigating solutions to various optimal control problems.}

Nowadays, there are many different definitions of fractional-order integrals 
or derivatives \cite{teo} and, in some sense, it is possible to consider a 
single broad class of fractional-order operators that include existing ones 
as particular cases \cite{aran}. This is important in analyzing a single operator 
rather than focusing on each individual one separately. Also, it enriches the 
subject of Pontryagin's maximum principle to handle a more general and wide class 
of fractional-order operators. For instance, a maximum principle is obtained 
for a combined fractional operator with a general analytic kernel in \cite{fays}. 
Also, some recent results for the Pontryagin maximum principle are investigated 
for fractional stochastic delayed systems with non-instantaneous 
impulses \cite{kumar}, for a degenerate fractional differential \mbox{equation \cite{banda},} 
and for distributed-order fractional derivatives \cite{faysn}. In contrast, 
the main aim of our current study is to utilize the idea of multi-order or 
incommensurate orders of derivatives in the definition of optimal control 
problems and then analyze their solutions. In doing so, 
we start with the most popular fractional model, which is the Caputo.

The structure of the article is as follows. In Section~\ref{sec:2:P}, we present some basic definitions 
and properties from fractional calculus. In Section~\ref{sectionR}, our contribution is given: 
we start by introducing the incommensurate non-local fractional optimal~control~problem, 
then we prove the continuity of solutions (Lemma~\ref{conti}), differentiability 
of perturbed trajectories (Lemma~\ref{differenti}), and, finally, the proof 
of the Pontryagin maximum principle (Theorem~\ref{theo}). In Section~\ref{sec:ex}, 
we give an example illustrating the application of the new Pontryagin maximum principle. 
We end with Section~\ref{sec:concl}, summarizing our work
and giving some perspectives for future work.


\section{Preliminaries}
\label{sec:2:P}

In this section, we briefly recall the necessary notions 
and results from fractional calculus. For more on the subject, 
we refer the interested readers to the books
\cite{Baleanu,MR1347689,Zhou}.

Let $(\alpha_i)_{i=1, \ldots, n} \in(0,1)$ be a multi-order 
of real numbers. In the sequel, we use the following notation: 
\[
L^{(\alpha)}\left([a, b], \mathbb{R}^n \right)
:=\left\{ x\in L^{1}\left([a, b], \mathbb{R}^n\right)
: I^{\alpha_i}_{a^{+}} x_i, I^{\alpha_i}_{b^{-}} x_i 
\in AC\left([a, b], \mathbb{R} \right)\right\},
\]
where $I^{\alpha_i}_{a^{+}}$ and $I^{\alpha_i}_{b^{-}}$ represent, respectively, 
the left and right Riemann--Liouville integrals of order $\alpha_i$.
We also use the notation $AC^{(\alpha)}\left([a, b], \mathbb{R}^n \right)$ to
represent the set of absolutely continuous functions that can be represented as 
\[
x_i(t)= x_i(a) + I^{\alpha_i}_{a^{+}}f(t) 
\, \text{ and }\, x_i(t)= x_i(b) + I^{\alpha_i}_{b^{-}} f(t),
\]
for some functions $f\in L^{\alpha_i}:=\left\{ x\in L^{1}\left([a, b], \mathbb{R}\right)
: I^{\alpha_i}_{a^{+}} x_i, I^{\alpha_i}_{b^{-}} x_i 
\in AC\left([a, b], \mathbb{R} \right) \right\}$.

\begin{Definition}
The left- and right-sided Riemann--Liouville incommensurate fractional derivative 
of orders $(\alpha_i)_{i=1, \ldots, n} \in(0,1)$ of a function 
$x\in L^{(\alpha)}$ are defined, respectively, by 
\[
D^{(\alpha)}_{a^{+}}x(t) 
= \begin{cases}
D^{\alpha_1}_{a^{+}}x_1(t),\\
D^{\alpha_2}_{a^{+}}x_2(t),\\
\quad \vdots  \quad \vdots\\
D^{\alpha_n}_{a^{+}}x_n(t),\\
\end{cases}; \qquad
D^{(\alpha)}_{b^{-}}x(t) 
= \begin{cases}
D^{\alpha_1}_{b^{-}}x_1(t),\\
D^{\alpha_2}_{b^{-}}x_2(t),\\
\quad \vdots  \quad \vdots\\
D^{\alpha_n}_{b^{-}}x_n(t),\\
\end{cases}
\]
with
\[
D^{\alpha_i}_{a^{+}}x_i(t)=\frac{d}{dt}\left( I^{1-\alpha_i}_{a^{+}} 
x_i(t)\right), \, \, D^{\alpha_i}_{b^{-}}x_i(t)
=-\frac{d}{dt}\left( I^{1-\alpha_i}_{b^{-}} x_i(t)\right),
\] 
where $I^{1-\alpha_i}_{a^{+}}$ and $I^{1-\alpha_i}_{b^{-}}$ represent, 
respectively, the left- and right-sided Riemann--Liouville 
fractional integrals of order $1-\alpha_i$.
\end{Definition}

\begin{Definition}
The left- and right-sided Caputo incommensurate fractional 
derivatives of order $(\alpha_i)_{i=1, \ldots, n} \in(0,1)$ 
of a function $x\in AC^{(\alpha)}$ are defined, respectively, by 
\[
{}^{\textsc c}D^{(\alpha)}_{a^{+}}x(t) 
= \begin{cases}
{}^{\textsc c}D^{\alpha_1}_{a^{+}}x_1(t),\\
{}^{\textsc c}D^{\alpha_2}_{a^{+}}x_2(t),\\
\quad \vdots  \quad \vdots\\
{}^{\textsc c}D^{\alpha_n}_{a^{+}}x_n(t),\\
\end{cases}; \qquad
{}^{\textsc c}D^{(\alpha)}_{b^{-}}x(t) 
= \begin{cases}
{}^{\textsc c}D^{\alpha_1}_{b^{-}}x_1(t),\\
{}^{\textsc c}D^{\alpha_2}_{b^{-}}x_2(t),\\
\quad \vdots  \quad \vdots\\
{}^{\textsc c}D^{\alpha_n}_{b^{-}}x_n(t),\\
\end{cases}
\]
where
\[
D^{\alpha_i}_{a^{+}}x_i(t)=  I^{1-\alpha_i}_{a^{+}}\left( 
\frac{d}{dt} x_i(t)\right), \quad D^{\alpha_i}_{b^{-}}x_i(t)
=-I^{1-\alpha_i}_{b^{-}}\left(  \frac{d}{dt}x_i(t)\right).
\] 
\end{Definition}

Note that integration by parts is a powerful tool when two functions are multiplied together, 
being useful for our purposes in the proof of the Pontryagin maximum principle.
In the sequel, the dot $\cdot$ is used for indicating scalar products.

\begin{Lemma}[Integration by parts formula \cite{MR1347689}]
\label{Lem:IbyP}
Let $x\in L^{(\alpha)}$ and $y\in AC^{(\alpha)}$. Then, 
\[
\int^b_a x(t)\cdot {}^{\textsc c}D^{(\alpha)}_{a^{+}}y(t)dt 
= \left[ y(t)\cdot I^{1-(\alpha)}_{b^{-}}x(t) \right]^b_a 
+ \int^b_a y(t)\cdot D^{(\alpha)}_{b^{-}}x (t)dt,
\]
where \[
I^{1-(\alpha)}_{b^{-}}x(t) 
= \begin{cases}
I^{1-\alpha_1}_{b^{-}}x_1(t),\\
I^{1-\alpha_2}_{b^{-}}x_2(t),\\
\quad \vdots  \quad \vdots\\
I^{1-\alpha_n}_{b^{-}}x_n(t).\\
\end{cases}
\]
\end{Lemma}
In what follows, we recall a generalized Gronwall inequality that 
is useful to prove the continuity and differentiability of perturbed trajectories.

\begin{Lemma}[Generalized Gronwall inequality \cite{gronwall}]
\label{lem:gronwall}
Let $\alpha$ be a positive real number and let $p(\cdot)$, 
$q(\cdot)$, and $u(\cdot)$ be non-negative continuous functions 
on $[a, b]$ with $q(\cdot)$ monotonic increasing on $[a, b)$.
If
\[
u(t)\leq p(t) + q(t)\int^t_a(t-s)^{\alpha-1}u(s)ds,
\]
then 
\[
u(t)\leq p(t) + \int^t_a \left[ \sum^{\infty}_{n=1}
\frac{\left(q(t)\Gamma(\alpha ) \right)^n}{\Gamma(n\alpha)} 
(t-s)^{n\alpha-1}p(s)\right]ds
\]
for all $t\in [a,b)$.
\end{Lemma}


\section{Main Results}
\label{sectionR}

In this section, our main concern is to find a control function 
$u \in L^{\infty}\left([a, b], \mathbb{R}^m \right)$ 
and its corresponding state trajectory 
$x\in AC^{(\alpha)}\left([a, b], \mathbb{R}^n \right)$, 
solution to the following incommensurate non-local 
fractional optimal~control~problem:
\begin{equation}
\label{pmp}
\begin{gathered}
J[x(\cdot), u(\cdot)]= \varphi(b, x(b)) 
+ \int^{b}_{a} L\left(t, x(t), u(t)\right)dt 
\longrightarrow \max,\\
{}^{\textsc c}D^{(\alpha)}_{{a}^{+}}x(t)
= f\left(t, x(t), u(t)\right), \quad t\in [a, b] \ a.e.,\\
x(\cdot) \in AC^{(\alpha)}, \quad u(\cdot) \in L^{\infty},\\
x(a)= x_a\in \mathbb{R}^n, \quad u(t)\in \Omega,
\end{gathered}
\end{equation}
where $\Omega$ is a closed subset of $\mathbb{R}^m$. 
The data functions $L: [a, b]\times \mathbb{R}^n \times \mathbb{R}^m 
\rightarrow \mathbb{R}$, $\varphi: \mathbb{R}^n \rightarrow \mathbb{R}$  
and $f: [a, b]\times \mathbb{R}^n \times \mathbb{R}^m \rightarrow \mathbb{R}^n$ 
are subject to the following assumptions:
\begin{itemize}
\item The function $\varphi$ is of class $\mathcal{C}^1$.
\item Functions $L$ and $f$ are continuous in all 
its three arguments and of class $\mathcal{C}^1$ with respect 
to  $x$ and, in particular, locally Lipschitz-continuous, that is, for every compact 
$B\subset \mathbb{R}^n$ and for all $x, y\in B$ there is $K> 0$ such that
$\left|L(t, x, u)-L(t, y, u) \right| \leq K \left\| x-y \right\|$ and
$ \left\| f(t, x, u)-f(t, y, u) \right\| \leq K \left\| x-y \right\|$.
\item There exists also $N>0$, such that 
$\displaystyle{\left\|\frac{\partial L(t, x, u)}{\partial x} \right\| \leq N}$ and
$\displaystyle{\left\| \frac{\partial f(t, x, u)}{\partial x} \right\| \leq N}$.

\item With respect to the control $u$, there exists $M>0$ such that 
\[
\left| L(t, x, u)\right| \leq M, \quad 
\parallel f(t, x, u)\parallel \leq M,
\quad \forall (t, x)\in [a, b]\times \mathbb{R}^n.
\]
\end{itemize}


\subsection{Needle-like Perturbation of the Optimal Control} 
\label{sub:sec:NC}

We will prove a first-order necessary optimality 
condition of the Pontryagin type, which is only sufficient 
in very particular cases \cite{Nikolskii}. Any proof of a necessary 
optimality condition, e.g., Fermat theorem about stationary points 
or the classical Pontryagin maximum \mbox{principle \cite{MR0186436},} 
begins by assuming the existence of a solution. We do the same here:
we assume $u^{*}(t)\in \Omega$ to be an optimal control to problem \eqref{pmp} 
for $t\in [a, b]$. The reader who is interested in the question of the existence of $u^{*}$
is referred to the recent paper \cite{MR4449935}.
Our aim, in this section, is to derive continuity and  differentiability 
properties of perturbed trajectories, which are crucial to proving the necessary 
optimality condition for the optimal control problem \eqref{pmp}. One way 
of achieving this is to perturb the optimal control by a needle-like variation 
and study the behavior of the corresponding state with respect to the optimal curve. 

Denote by $\mathcal{L}\left[ F(\cdot)\right]$ the set of all Lebesgue 
points in $[a, b)$ of the essentially bounded functions 
$t\mapsto f(t, x(t), u(t))$ and $t\mapsto L(t, x(t), u(t))$. Then, 
for $(\tau, v) \in \mathcal{L}\left[ F(\cdot)\right] \times \Omega$, 
and for every $\theta \in [0, b-\tau)$, let us consider the needle-like 
variation $u^{\theta} \in L^{\infty}\left([a, b], \mathbb{R}^n \right)$ 
associated to the optimal control $u^{*}$, which is given by
\begin{equation}
\label{needleVar}
u^{\theta}(t)= 
\begin{cases}
u^{*}(t) \quad \text{ if } \quad t\not\in [\tau,\tau+\theta),\\
v \qquad \quad \text{ if } \quad t\in [\tau, \tau + \theta),
\end{cases}
\end{equation}
for almost every $t\in [a, b]$.

\begin{Lemma}[Continuity of solutions]
\label{conti} 
For any $(\tau, v) \in \mathcal{L}\left[ F(\cdot)\right] \times \Omega$, 
denote by $x^{\theta}$ the corresponding state trajectory to the needle-like 
variation $u^{\theta}$, that is, the state solution of 
\begin{equation}\label{continuity}
{}^{\textsc c}D^{(\alpha)}_{a+} x^{\theta}(t)
= f\left(t, x^{\theta}(t), u^{\theta}(t)\right), \quad x^{\theta}(a)= x_a.
\end{equation}

Then, the state $x^{\theta}$ converges uniformly to the optimal state 
trajectory $x^{*}$ whenever $\theta$ tends to zero.
\end{Lemma}

\begin{proof}
By definition of incommensurate Caputo derivative, we have
\[
{}^{\textsc c}D^{(\alpha)}_{a+}\left( x^{\theta}(t) 
- x^{*}(t)\right) 
=  f(t, x^{\theta}(t), u^{\theta}(t)) 
- f\left(t, x^{*}(t), u^{*}(t)\right).
\]

Then, the integral representation is obtained as 
\[
x^{\theta}(t) 
- x^{*}(t)
= \int_0^t {\rm diag}\left(
\frac{(t-s)^{\alpha_i - 1}}{\Gamma(\alpha_i)} \right)
\cdot \Big( f(t, x^{\theta}(t), u^{\theta}(t)) 
- f\left(t, x^{*}(t), u^{*}(t)\right)\Big)ds.
\]

Moreover, the function $\displaystyle{ \alpha_{i} \mapsto 
\frac{c^{\alpha_{i}-1}}{\Gamma(\alpha_{i})} }$ is continuous 
on $\displaystyle{[\underset {1\leq i \leq n} \min \alpha_i, 1]}$, 
where $c$ is a non-zero constant. Thus, by the extreme value theorem 
due to Weierstrass, it attains a maximum. Hence, there exists 
$\bar \alpha \in [\underset {1\leq i \leq n} \min \alpha_i, 1]$ 
such that $ \displaystyle{\frac{c^{\alpha_{i}-1}}{\Gamma(\alpha_{i})} 
\leq \frac{c^{\bar \alpha-1}}{\Gamma(\bar \alpha)}}$. This leads to 
\begin{multline*}
x^{\theta}(t) 
- x^{*}(t)
\leq   \int_0^t 
\frac{(t-s)^{\bar \alpha - 1}}{\Gamma(\bar \alpha)}
\cdot \Big( f(t, x^{\theta}(t), u^{\theta}(t)) 
- f\left(t, x^{*}(t), u^{*}(t)\right)\Big)ds\\
= I^{\bar \alpha}_{a+} \Big( f(t, x^{\theta}(t), u^{\theta}(t)) 
- f\left(t, x^{*}(t), u^{*}(t)\right)\Big).
\end{multline*}

Further, the following relation holds for  function $f$, 
\begin{multline*}
f(t, x^{\theta}(t), u^{\theta}(t)) 
- f\left(t, x^{*}(t), u^{*}(t)\right)
= \{ f(t, x^{\theta}(t), u^{\theta}(t)) - f(t, x^{*}(t), u^{\theta}(t))\} \\
+ \{f(t, x^{*}(t), u^{\theta}(t)) - f\left(t, x^{*}(t), u^{*}(t)\right)\}.
\end{multline*}

Using the triangular inequality, and noticing that  $u^{\theta}$ and $u^{*}$ are different only 
on $[\tau, \tau+\theta)$, we obtain
\begin{multline*}
\parallel x^{\theta}(t)-x^{*}(t) \parallel \leq 
I^{\bar{\alpha}}_{a^{+}}\left( \parallel f(t, x^{\theta}(t), u^{\theta}(t)) 
-f(t, x^{*}(t), u^{\theta}(t)) \parallel \right)\\
+  I^{\bar{\alpha}}_{\tau^{+}}\left( 
\parallel f(t, x^{*}(t), u^{\theta}(t))
- f\left(t, x^{*}(t), u^{*}(t)\right) \parallel \right).
\end{multline*}

Next, from data assumptions, $K$ is a Lipschitz constant of $f$ and $M$ 
is an upper bound of $f$ with respect to the control function. 
Thus, it follows that
\[
\parallel x^{\theta}(t)-x^{*}(t) \parallel \leq  
K I^{\bar{\alpha}}_{a^{+}}\left( \parallel x^{\theta}(t)
-x^{*}(t) \parallel \right) + 
M\frac{\theta^{\bar \alpha}}{\Gamma(\bar \alpha + 1)}.
\]

Now, by applying the generalized Gronwall inequality
of Lemma~\ref{lem:gronwall}, it follows that
\begin{multline*}
\parallel x^{\theta}(t)-x^{*}(t) \parallel 
\leq \frac{M\theta^{\bar{\alpha}}}{\Gamma(\alpha + 1)}\left[   
1+ \int_a^t \sum^{\infty}_{n=1}\frac{K^n}{\Gamma(n\bar{\alpha})} 
(t-s)^{n\bar{\alpha}-1}ds \right]\\
 \leq \frac{M \theta^{\bar \alpha}}{\Gamma(\alpha + 1)} E_{\alpha, 1}\left( K(b-a)^{\alpha}\right),
\end{multline*}
where  $E_{\alpha, 1}$ is the Mittag--Leffler function of parameter $\bar{\alpha}$ \cite{MR4650052}.
Hence, we obtain that $\left\| x^{\theta}-x^{*} \right\|$ converges to zero on $[a, b]$, which ends the proof.
\end{proof}

\begin{Lemma}[Differentiability of the perturbed trajectory]
\label{differenti}
Suppose that $(x^{*}, u^{*})$ is an optimal pair to problem \eqref{pmp}. 
Then, for $(\tau, v) \in \mathcal{L}\left[ F(\cdot)\right] \times \Omega$, 
the quotient variational trajectory 
$\displaystyle{\frac{x^{\theta}(\cdot)-x^{*}(\cdot)}{\theta}}$ converges 
uniformly on $[\tau + \theta, b]$ to $\eta(\cdot)$ when $\theta$ tends 
to zero, where $\eta(\cdot)$ is the unique solution 
to the incommensurate left Caputo fractional Cauchy problem
\begin{equation}
\label{eqLinear}
\begin{cases}
\displaystyle{ {}^{\textsc c}D^{(\alpha)}_{\tau +} \eta(t)
=  \frac{\partial f(t, x^{*}(t), u^{*}(t))}{\partial x} \cdot \eta(t)}, 
\quad t\in ]\tau, b],\\[3mm]
\displaystyle{I^{1-\bar{\alpha}}_{{\tau}^{+}} \eta(\tau)
= \left(f(\tau, x^{*}(\tau), v)
-f(\tau, x^{*}(\tau), u^{*}(\tau))\right)},
\end{cases}
\end{equation}
with $\displaystyle{\bar \alpha \in [\underset {1\leq i \leq n} \min \alpha_i, 1]}$, 
such that $\displaystyle{\frac{c^{\alpha_i -1}}{\Gamma(\alpha_i)}
\leq \frac{c^{\bar \alpha -1}}{\Gamma(\bar \alpha)}, \, 
\text{ for some non-zero constant } c}$ and $i=1,\ldots, n$.
\end{Lemma} 

\begin{proof}
Set $\displaystyle{z^{\theta}(t)= \frac{x^{\theta}(t)-x^{*}(t)}{\theta}-\eta(t)}$ 
for all $t\in [\tau +\theta, b]$. Our aim is to prove that $z^{\theta}$ 
converges uniformly to zero whenever $\theta \rightarrow 0$. 

The integral representation of $z^{\theta}$ is obtained as follows:
\begin{multline*}
z^{\theta}(t)= \int^{\tau +\theta}_{\tau} {\rm diag}\left(
\frac{(t-s)^{\alpha_i - 1}}{\Gamma(\alpha_i)} \right)
\cdot\frac{f(s, x^{\theta}(s), v)-f(s, x^{*}(s), v)}{\theta}ds\\
+ \int^{\tau +\theta}_{\tau} {\rm diag}\left(
\frac{(t-s)^{\alpha_i - 1}}{\Gamma(\alpha_i)} \right)
\cdot\frac{f(s, x^{*}(s), v)-f(s, x^{*}(s), u^{*}(s))}{\theta}ds\\
-\frac{1}{\Gamma(\bar \alpha)}(t-\tau)^{\bar{\alpha}-1}\left( 
f(\tau, x^{*}(\tau), v)-f(\tau, x^{*}(\tau), u^{*}(\tau))\right)\\
+\int^t_{\tau +\theta} {\rm diag}\left(
\frac{(t-s)^{\alpha_i - 1}}{\Gamma(\alpha_i)} \right)\cdot\left(
\frac{f(s, x^{\theta}(s), u^{*}(s))-f(s, x^{*}(s), u^{*}(s))}{\theta} \right.\\
\left. -\frac{\partial f(s, x^{*}(s), u^{*}(s))}{\partial x}
\times \frac{x^{\theta}(s)-x^{*}(s)}{\theta}\right) ds\\
\int^{\tau +\theta}_{\tau} {\rm diag}\left(
\frac{(t-s)^{\alpha_i - 1}}{\Gamma(\alpha_i)} \right)
\cdot \frac{\partial f(s, x^{*}(s), u^{*}(s))}{\partial x}
\times \eta(s)ds\\
+\int^t_{\tau +\theta} {\rm diag}\left(
\frac{(t-s)^{\alpha_i - 1}}{\Gamma(\alpha_i)} \right)
\cdot \frac{\partial f(s, x^{*}(s), u^{*}(s))}{\partial x}
\times z^{\theta}(s)ds.
\end{multline*}

Note that by the existence property of $\bar \alpha$, we have 
that $\displaystyle{\frac{(t-s)^{\alpha_i -1}}{\Gamma(\alpha_i)}
\leq \frac{(t-s)^{\bar \alpha -1}}{\Gamma(\bar \alpha)}}$ 
for all $i=1,\ldots,n$. Thus, we can deduce that
\begin{multline}
\label{eq:Intrepres}
z^{\theta}(t) \leq \frac{1}{\Gamma(\bar \alpha)}\int^{\tau 
+\theta}_{\tau} (t-s)^{\bar \alpha - 1}
\frac{f(s, x^{\theta}(s), v)-f(s, x^{*}(s), v)}{\theta}ds\\
+ \frac{1}{\Gamma(\bar \alpha)}\int^{\tau +\theta}_{\tau} 
(t-s)^{\bar \alpha - 1}\frac{f(s, x^{*}(s), v)-f(s, x^{*}(s), u^{*}(s))}{\theta}ds\\
-\frac{1}{\Gamma(\bar \alpha)}(t-\tau)^{\bar{\alpha}-1}\left( 
f(\tau, x^{*}(\tau), v)-f(\tau, x^{*}(\tau), u^{*}(\tau))\right)\\
+\frac{1}{\Gamma(\bar \alpha)}\int^{t}_{\tau +\theta} (t-s)^{\bar \alpha - 1}\left(
\frac{f(s, x^{\theta}(s), u^{*}(s))-f(s, x^{*}(s), u^{*}(s))}{\theta} \right.\\
\left. -\frac{\partial f(s, x^{*}(s), u^{*}(s))}{\partial x}
\times \frac{x^{\theta}(s)-x^{*}(s)}{\theta}\right) ds\\
-\frac{1}{\Gamma(\bar \alpha)}\int^{\tau +\theta}_{\tau} 
(t-s)^{\bar \alpha - 1} \frac{\partial f(s, x^{*}(s), u^{*}(s))}{\partial x}
\times \eta(s)ds\\
+\frac{1}{\Gamma(\bar \alpha)}\int^{t}_{\tau +\theta} 
(t-s)^{\bar \alpha - 1} \frac{\partial f(s, x^{*}(s), u^{*}(s))}{\partial x}
\times z^{\theta}(s)ds.
\end{multline}

In line with the work of \cite{bourdin}, precisely the proof of 
Proposition 3.3 of this reference,  we obtain that each term appearing 
in the right hand side of \eqref{eq:Intrepres} 
is bounded, which yields the following estimate:
\[
\lVert z^{\theta}(t)\rVert \leq \Theta^{\theta}_1(t-(\tau +\theta))^{\bar \alpha -1}
E_{\bar \alpha, \bar \alpha^{'}}\left( N(b-a)^{\bar \alpha}\right) 
+ \Theta^{\theta}_2(t-(\tau +\theta))^{\bar \alpha^{'} -1}
E_{\bar \alpha, \bar \alpha^{'}}\left( N(b-a)^{\bar \alpha}\right),
\]
where functions $\Theta^{\theta}_1$ and $\Theta^{\theta}_2$ both converge 
uniformly to zero whenever $\theta$ tends to zero. Hence, we conclude that 
$z^{\theta}(t)$ converges uniformly to zero as $\theta$ 
goes to zero, which is the desired result. 
\end{proof}


\subsection{Pontryagin's Maximum Principle for Problem \eqref{pmp}}
\label{subsec:PMP}

The fractional Pontryagin maximum principle has many applications, 
for example in Engineering, Economics and Health \cite{MR4537621,MR4625042}.
Here, we state and prove the main result of our work: a Pontryagin maximum principle
for the incommensurate fractional-order optimal control problem \eqref{pmp}.

\begin{Theorem}[Pontryagin Maximum Principle for \eqref{pmp}]
\label{theo}
If $(x^{*}(\cdot), u^{*}(\cdot))$ is an optimal pair for \eqref{pmp}, 
then there exists $\lambda \in L^{(\alpha)}$, called the adjoint function variable, 
such that the following conditions hold for all $t$ in the interval $[a, b]$:
\begin{itemize}
\item the maximality condition
\begin{equation}
\label{maxcondition}
H(t, x^{*}(t), u^{*}(t), \lambda(t))
= \underset{\omega \in \Omega} \max \, H\left(t, x^{*}(t), \omega, \lambda(t)\right);
\end{equation}
\item the adjoint system
\begin{equation}
\label{adj}
 D^{(\alpha)}_{b-} \lambda(t)
= \frac{\partial H}{\partial x}(t,x^{*}(t), u^{*}(t), \lambda(t));
\end{equation}
\item the transversality condition
\begin{equation}
\label{trans}
I^{1-(\alpha)}_{b^{-}}\lambda(b)=\frac{\partial \varphi}{\partial x}(b, x^{*}(b)),
\end{equation}
\end{itemize}
where the Hamiltonian function $H$ is defined by
\begin{equation}
\label{eq:Hamilt}
H(t, x, u, \lambda)= L(t, x, u) + \lambda \cdot f(t, x, u).
\end{equation}
\end{Theorem}

\begin{proof}
Let $x^{\theta}(t)$ be the corresponding state trajectory to the needle-like 
variation $u^\theta$ defined by \eqref{needleVar}. Observe that, 
by integration by parts, we have for $\lambda(\cdot)\in L^{(\alpha)}$, 
\begin{equation}
\int^b_a \lambda(t)\cdot {}^{\textsc c}D^{(\alpha)}_{a^{+}}x^{\theta}(t)dt 
= \left[x^{\theta}(t)\cdot I^{1-(\alpha)}_{b-}\lambda(t) \right]^b_a 
+ \int^b_a x^{\theta}(t)\cdot D^{(\alpha)}_{b-}\lambda(t)dt.
\end{equation}

This relation can be added to the objective functional 
at $(x^{\theta}, u^{\theta})$ defined by 
$$
\displaystyle{ J(x^{\theta}, u^{\theta}) = \varphi(b, x^{\theta}(b))
+ \int_a^b L\left( t, x^{\theta}(t), u^{\theta}(t)\right)ds},
$$ 
meaning that
\vspace{-6pt}
\begin{adjustwidth}{-\extralength}{0cm}
\begin{multline}
J(x^{\theta}, u^{\theta})= \varphi(b, x^{\theta}(b))
+ \int_a^b \Big[L\left( t, x^{\theta}(t), u^{\theta}(t)\right) 
+ \lambda(t)\cdot {}^{\textsc c}D^{(\alpha)}_{a^{+}}x^{\theta}(t)
- x^{\theta}(t)\cdot D^{(\alpha)}_{b-}\lambda(t)\Big]dt\\
x^{\theta}(b)\cdot I^{1-(\alpha)}_{b-}\lambda(b) 
+ x^{\theta}(a)\cdot I^{1-(\alpha)}_{a}\lambda(a),
\end{multline}
\end{adjustwidth}
which, by substituting \eqref{continuity} to this latter expression, leads to
\begin{multline*}
J(x^{\theta}, u^{\theta}) = \varphi(b, x^{\theta}(b))
+ \int_a^b \left[ H\left( t, x^{\theta}(t), u^{\theta}(t), 
\lambda (t) \right)- x^{\theta}(t)\cdot D^{(\alpha)}_{b-}\lambda(t)\right]dt\\ 
- x^{\theta}(b)\cdot I^{1-(\alpha)}_{b-}\lambda(b) 
+ x_a\cdot I^{1-(\alpha)}_{a}\lambda(a),
\end{multline*}
where $H(t, x, u, \lambda)= L(t, x, u) + \lambda \cdot f(t, x, u)$.
Next, we write down the Taylor expansion 
\vspace{-6pt}
\begin{adjustwidth}{-\extralength}{0cm}
\begin{multline*}
\varphi(b, x^{\theta}(b)) = \varphi(b, x^{*}(b)) 
+ \left( x^{\theta}(b) - x^{*}(b)\right)
\cdot \frac{\partial \varphi}{\partial x}(b, x^{*}(b)) 
+ \circ\left(\parallel x^{\theta}-x^{*}\parallel \right);\\
H\left( t, x^{\theta}(t), u^{\theta}(t), \lambda(t) \right) 
=  H\left( t, x^{*}(t), u^{\theta}(t), \lambda(t) \right)
+ \left( x^{\theta}(t) - x^{*}(t)\right)\cdot \frac{\partial H}{\partial x}
\left( t, x^{*}(t), u^{\theta}(t), \lambda(t)\right)\\ 
+ o\left(\parallel x^{\theta}-x^{*}\parallel \right).
\end{multline*}
\end{adjustwidth}

Note that, by the continuity Lemma~\ref{conti}, we have the uniform 
convergence of  $\|x^{\theta}-x^{*}\| \rightarrow 0$ whenever 
$\theta \rightarrow 0$. Thus, the residue term in the Taylor 
expansion can be expressed as a function of $\theta$. Therefore, 
we can evaluate the quotient $\displaystyle{ \frac{J(x^{\theta}, 
u^{\theta}) -J(x^{*}, u^{*})}{\theta}  := \delta J}$ as follows:
\vspace{-12pt}
\begin{adjustwidth}{-\extralength}{0cm}
\begin{multline*}
\delta J = \frac{x^{\theta}(b)- x^{*}(b)}{\theta}
\cdot \frac{\partial \varphi}{\partial x}(b, x^{*}(b)) 
+  \int_a^b \frac{H\left( t, x^{*}(t), u^{\theta}(t), \lambda(t) \right) 
- H\left( t, x^{*}(t), u^{*}(t), \lambda(t) \right)}{\theta}dt \\
+ \int^b_a \Big(  \frac{\partial H}{\partial x}\left( t, x^{*}(t), 
u^{\theta}(t), \lambda(t)\right)- D^{(\alpha)}_{b-}\lambda(t) \Big)
\cdot \frac{ x^{\theta}(t) - x^{*}(t)}{\theta}dt - \frac{x^{\theta}(b)
- x^{*}(b)}{\theta}\cdot I^{1-(\alpha)}_{b-}\lambda(b)\\+ o(\theta)(1 + b-a).
\end{multline*}
\end{adjustwidth}

Now, by the differentiability Lemma~\ref{differenti}, we have that
$\frac{ x^{\theta}(t) - x^{*}(t)}{\theta}$ converges uniformly 
to $\eta(t)$ when $\theta$ tends to zero. Therefore, the limit 
of $\delta J$ when $\theta$ tends to zero can be expressed as 
\begin{multline*}
\lim_{\theta \rightarrow 0} \delta J = \eta(b)\cdot 
\frac{\partial \varphi}{\partial x}(b, x^{*}(b)) 
+ \lim_{\theta \rightarrow 0}\int_a^b  \frac{H\left( t, x^{*}(t), 
u^{\theta}(t), \lambda(t) \right) - H\left( t, x^{*}(t), u^{*}(t), 
\lambda(t) \right)}{\theta}dt \\
+ \int^b_a \Big(  \frac{\partial H}{\partial x}\left( t, x^{*}(t), u^{*}(t), 
\lambda(t) \right)- D^{(\alpha)}_{b-}\lambda(t) \Big)
\cdot \eta(t)dt - \eta(b) \cdot I^{1-(\alpha)}_{b-}\lambda(b).
\end{multline*}

Next, we fix
\[
D^{(\alpha)}_{b-}\lambda(t) = \frac{\partial H}{\partial x}
\left( t, x^{*}(t), u^{*}(t), \lambda(t) \right) \quad \textup{ with } 
\quad I^{1-(\alpha)}_{b-}\lambda(b) =\frac{\partial \varphi}{\partial x}(b, x^{*}(b)),
\]
that is, the adjoint Equation \eqref{adj} and the transversality 
condition \eqref{trans}. Thus, we are left with 
\[
\lim_{\theta \rightarrow 0} \delta J =  \lim_{\theta \rightarrow 0}
\int_a^b  \frac{H\left( t, x^{*}(t), u^{\theta}(t), \lambda(t) \right) 
- H\left( t, x^{*}(t), u^{*}(t), \lambda(t) \right)}{\theta}dt.
\]

Moreover, recalling that
$
u^{\theta}(t)= 
\begin{cases}
u^{*}(t) \quad \text{ if } \quad t\not\in [\tau, \tau+\theta);\\
v \qquad \quad \text{ if } \quad t\in [\tau, \tau +\theta),
\end{cases}
$
for almost every $t\in [a, b]$, it follows that
\[
\lim_{\theta \rightarrow 0} \delta J = \lim_{\theta 
\rightarrow 0^{+}}\frac{1}{\theta}\int^{\tau+\theta}_{\tau} \left[ 
H(s, x^{*}(s), v, \lambda(s))-H(s, x^{*}(s), u ^{*}(s), \lambda(s)) \right]ds .
\]

However, notice that $\tau$ is a Lebesgue point of 
$$
H(s, x^{*}(s), v, \lambda(s))-H(s, x^{*}(s), u ^{*}(s), \lambda(s)) := \psi(s).
$$

Thus, from the Lebesgue differentiation property,
\[
\left| \frac{1}{\theta}\int^{\tau +\theta}_{\tau}\psi(s)ds -\psi(\tau)\right| 
= \left| \frac{1}{\theta}\int^{\tau +\theta}_{\tau }\left(\psi(s)-\psi(\tau)\right)ds \right| 
\leq \frac{1}{\theta}\int^{\tau +\theta}_{\tau}\left|\psi(s)-\psi(\tau)\right|ds,
\]
and we have that the right-hand side tends to zero for almost every point $\tau$. 
As a consequence,
\begin{equation}
\begin{split}
\lim_{\theta \rightarrow 0^{+}}\frac{1}{\theta}\int^{\tau +\theta}_{\tau} \left[ 
H(s, x^{*}(s), v, \lambda(s))-H(s, x^{*}(s), u ^{*}(s), \lambda(s)) \right]ds\\
=H(\tau, x^{*}(\tau), v, \lambda(\tau))-H(\tau, x^{*}(\tau), u ^{*}(\tau), \lambda(\tau)).
\end{split}
\end{equation}

Further, by the optimality assumption of the pair 
\[ 
(x^{*}, u^{*}), \text{ one has also } \lim_{\theta \rightarrow 0} \delta J \leq 0,
\]
altogether, we obtain that
\[
H(\tau, x^{*}(\tau), v, \lambda(\tau))-H(\tau, x^{*}(\tau), u ^{*}(\tau), \lambda(\tau))\leq 0.
\] 

Finally, because $\tau$ is an arbitrary Lebesgue point of the control $u^{*}$ 
and $v$ is an arbitrary element of the set $\Omega$, 
it follows that the relation 
\[
H(t, x^{*}(t), u^{*}(t), \lambda(t))
= \underset{\omega \in \Omega} 
\max \, H\left(t, x^{*}(t), \omega, \lambda(t)\right)
\]
holds at all Lebesgue points, which ends the proof.
\end{proof}

Theorem~\ref{theo} gives a necessary optimality condition 
that provides an algorithm that can be used to solve
general Bolza-type fractional optimal control problems 
that depend on multi-order Caputo fractional derivatives:
given a problem \eqref{pmp},
\begin{itemize}
\item[(i)] one writes the associated Hamiltonian \eqref{eq:Hamilt}; 
\item[(ii)] we use the maximality condition \eqref{maxcondition} 
to obtain an expression of the optimal controls in terms of the state and adjoint variables; 
\item[(iii)] we substitute the expressions obtained in step (ii) in the  adjoint system \eqref{adj};
\item[(iv)] finally, we solve the system obtained in step (iii) together with the initial conditions 
$x(a) = x_a$ and the transversality condition \eqref{trans}. 
\end{itemize}

A simple example of the usefulness of our result is given in Section~\ref{sec:ex}.


\section{An Illustrative Example}
\label{sec:ex}

Let us study the following optimal control multi-order fractional Caputo problem:
\begin{equation}
\label{exple}
\begin{gathered}
x_2(5) + \int^{5}_{1} \left(1+\exp(2u(t))\right)dt 
\longrightarrow \max,\\
x=(x_1, x_2) \in AC^{(\frac{1}{3}, \frac{1}{2})}\left([1, 5], 
\mathbb{R}^2 \right), \quad u  \in L^{\infty}\left([1, 5], \mathbb{R} \right),\\
\begin{cases}
{}^{\textsc c}D^{\frac{1}{3}}_{{1}^{+}}x_1(t)
= 1-\exp(2u(t)), \\[3mm]
{}^{\textsc c}D^{\frac{1}{2}}_{{1}^{+}}x_2(t)
= x_1(t), 
\end{cases}\\
x(1)= (1, 1), \quad u(t)\in [-2, 7].
\end{gathered}
\end{equation}

The Hamiltonian function and the running cost 
of this problem \eqref{exple} are given, respectively, by
\[
H(t, x, u, \lambda)= 1+\exp(2u) + \lambda_1(1-\exp(2u))
+ \lambda_2x_1, \, \, \text{ and } \varphi(b, x(b))= x_2(5).
\]

We do not know if the problem has a solution or not but, 
if the problem has a solution, then it must satisfy 
the necessary optimality condition is given by our Theorem~\ref{theo}. 
Precisely, if $(x^{*}, u^{*})$ is an optimal pair solution to problem \eqref{exple}, 
then by application of our Theorem~\ref{theo}, there exists an adjoint 
function $\lambda \in L^{(\frac{1}{3}, \frac{1}{2})}$, satisfying
\begin{equation*}
\begin{cases}
D^{\frac{1}{3}}_{5^{-}}\lambda_1(t)= \frac{\partial H}{\partial x_1}
= \lambda_2, \quad I^{1-\frac{1}{3}}_{5^{-}}=\frac{\partial \varphi}{\partial x_1} =0;\\
D^{\frac{1}{2}}_{5^{-}}\lambda_2(t)= \frac{\partial H}{\partial x_2}=0, 
\quad I^{1-\frac{1}{2}}_{5^{-}}=\frac{\partial \varphi}{\partial x_1} =1.
\end{cases}
\end{equation*}

We integrate, to obtain that
\[
\lambda_1(t)\frac{(5-t)^{\frac{2}{3}} - 1}{\Gamma(\frac{2}{3})}, 
\quad \text{ and } \quad \lambda_2(t)= \frac{(5-t)^{\frac{1}{2}-1}}{\Gamma(\frac{1}{2})}.
\]

Moreover, from the maximality condition \eqref{maxcondition}, it yields that
\[
u^{*}(t) \in \underset{v\in [-2, 7]}{\arg} \left\langle 
\begin{pmatrix}
1-\exp(2v(t))\\
1+\exp(2v(t))
\end{pmatrix} 
\cdot 
\begin{pmatrix}
\lambda_1\\
1
\end{pmatrix} \right\rangle,
\]
for almost $t\in [1, 5]$. Therefore, using the 
classical Cauchy--Schwartz inequality, we obtain
\[
\begin{pmatrix}
1-\exp(2u^{*}(t))\\
1+\exp(2u^{*}(t))
\end{pmatrix} = \frac{1}{|(\lambda_1(t), 1)|_{\mathbb{R}^2}}\begin{pmatrix}
\lambda_1\\
1
\end{pmatrix}.
\]

This leads to \[
\tanh\left(u^{*}(t) \right) = \frac{1}{\lambda_1(t)}=\Gamma\left(\frac{2}{3}\right)(5-t)^{\frac{1}{3}}.
\]

Finally, we obtain that $u^{*}$ is given by
\begin{equation}
\label{eq:extremal:control:ex}
u^{*}(t)= \arctan h \left( \Gamma\left(\frac{2}{3}\right)(5-t)^{\frac{1}{3}}\right), 
\text{ for almost } t\in [1, 5].
\end{equation}

Note that \eqref{eq:extremal:control:ex} 
is just a Pontryagin extremal (candidate).

For practical applications, the problems are difficult and one needs to use numerical methods.
Many software packages for computing fractional-order derivatives and solving fractional
differential systems are now available. We refer the interested reader to \cite{Garrappa}
and references therein.


\section{Conclusions}
\label{sec:concl}

In this paper, we have studied a general Bolza-type fractional optimal control problem 
that depends on multi-order Caputo fractional derivatives. We have established 
a Pontryagin maximum principle for this incommensurate fractional-order problem. 
Our approach starts with a sensitivity analysis from which we prove the continuity 
and differentiability of perturbed trajectories to the optimal state. 
An illustrative example shows the applicability of our main result (Theorem~\ref{theo}), 
which provides a Pontryagin maximum principle for incommensurate fractional-order
optimal control problems.

Recent applications of fractional mathematical modeling have shown the importance of considering
incommensurate orders in infectious disease dynamics \cite{faysk}. This might permit greater
flexibility in capturing the heterogeneous nature of disease dynamics, accounting for factors 
such as population demographics or social behaviors of \mbox{individuals \cite{jan}.} Further, it has an
advantage when introducing a control function in such models, by the fact that the more accurate the
model, the better the control. Therefore, for future work, it would be interesting to emphasize 
numerical methods for incommensurate-order problems in order to handle applications on existing
models of multi-order derivatives.


\vspace{6pt} 

\authorcontributions{Conceptualization, F.N. and D.F.M.T.; 
methodology, F.N. and D.F.M.T.; 
validation, F.N. and D.F.M.T.; 
formal analysis, F.N. and D.F.M.T.; 
investigation, F.N. and D.F.M.T.;  
writing---original draft preparation, F.N. and D.F.M.T.; 
writing---review and editing, F.N. and D.F.M.T. 
All authors have read and agreed to the published version of the manuscript.}

\funding{This research was funded by the Funda\c{c}\~{a}o para a Ci\^{e}ncia 
e a Tecnologia, I.P. (FCT, Funder ID = 50110000187) 
under Grants UIDB/04106/2020 and UIDP/04106/2020.}

\dataavailability{No new data were created or analyzed in this study. 
Data sharing is not applicable to this article.}

\acknowledgments{The authors are very grateful 
to three Reviewers for several comments, questions, 
and suggestions, that helped them to improve the submitted paper.}

\conflictsofinterest{The authors declare no conflicts of interest. 
The funders had no role in the design of the study; in the collection, 
analyses, or interpretation of data; in the writing of the manuscript; 
or in the decision to publish the results.} 


\begin{adjustwidth}{-\extralength}{0cm}
	
\reftitle{References}


\PublishersNote{}
\end{adjustwidth}

\end{document}